\documentclass{amsart}

\usepackage[inactive]{srcltx} 
\usepackage{amsmath, amsthm, amscd, amsfonts, amssymb, graphicx, color}
\usepackage[all]{xy}

\newtheorem{thm}{Theorem}[section]
\newtheorem{cor}[thm]{Corollary}
\newtheorem{lem}[thm]{Lemma}

\theoremstyle{definition}

\newtheorem{exm}[thm]{Example}
\theoremstyle{remark}
\newtheorem{rem}[thm]{Remark}

\begin{document}

\title[Lie triple maps on generalized matrix algebras]
{Lie triple maps on generalized matrix algebras}

\author{ B. Fadaee}
\thanks{{\scriptsize
\hskip -0.4 true cm \emph{MSC(2010)}: 16W25, 47B47, 15A78, , 16W10.
}}

\address{Department of
Mathematics, University of Kurdistan, P. O. Box 416, Sanandaj,
Iran.}

\email{behroozfadaee@yahoo.com; b.fadaee@sci.uok.ac.ir}

\address{}

\email{}

\thanks{}

\thanks{}



\date{}

\dedicatory{}

\commby{}


\begin{abstract}
In this article, we introduce the notion of Lie triple centralizer as follows. Let $\mathcal{A}$ be an algebra, and $\phi: \mathcal{A}\rightarrow \mathcal{A}$ be a linear mapping. We say that $\phi$ is a Lie triple centralizer whenever $\phi([[a,b],c])=[[\phi(a),b],c] $ for all $a,b,c\in \mathcal{A}$. Then we characterize the general form of Lie triple centralizers on a generalized matrix algebra $\mathcal{U}$ and under some mild conditions on $\mathcal{U}$ we present the necessary and sufficient conditions for a Lie triple centralizer to be proper. As an application of our results, we characterize generalized Lie triple derivations on generalized matrix algebras.
\end{abstract}

\maketitle
\section{Introduction}
Throughout the paper, all algebras and modules will be over the unital commutative ring $ \mathcal{R}$. Let $ \mathcal{A} $ be an algebra. A linear mapping $ \varphi : \mathcal{A} \to \mathcal{A} $ is called a \textit{Lie centralizer} if 
\[ \varphi ( [ a, b])= [ \varphi (a) , b ] \quad (a,b\in \mathcal{A}),\] 
where $[a, b]=ab -ba $ is the Lie brackets. It is easily checked that $ \varphi $ is a Lie centralizer on $ \mathcal{A} $ if and only if $\varphi ( [ a, b])= [ a, \varphi (b) ] $ for any $a , b \in \mathcal{A} $.The Lie centralizer on the Lie algebras is also called Lie centroid, which are important in studying the Lie structure of algebras (see \cite{jac, mc}). Recently, the study of the structure of Lie centralizers on algebras has been considered. In \cite{fo}, Fo\v{s}ner and Jing have described the non-additive Lie centralizers on triangular algebras. Jabeen in \cite{jab} has described Lie centralizers on generalized matrix algebras, and in \cite{liu2} non-linear Lie centralizers on generalized matrix algebra have been studied. In \cite{beh}, the authors have studied the characterization of Lie centralizers on non-unital triangular algebras through zero products, and in \cite{gh1} Lie centralizers at zero products on a class of operator algebras have been characterized. 
\par
Another important classes of mappings on algebras are derivations and Lie derivations, and their generalizations. Let $ \mathcal{A} $ be an algebra. A linear map $\delta:\mathcal{A} \rightarrow \mathcal{A}$ is said to be a \textit{derivation} if
\[\delta(ab)=\delta(a)b+a\delta(b)\quad (a,b\in \mathcal{A}).\]
A linear map $\sigma:\mathcal{A} \rightarrow \mathcal{A}$ is said to be a \textit{Lie derivation} if
\[\sigma([a,b])=[\sigma(a),b]+[a,\sigma(b)]\quad (a,b\in \mathcal{A}).\]
A linear map $\Delta:\mathcal{A} \rightarrow \mathcal{A}$ is said to be a \textit{generalized Lie derivation associated with the Lie derivation $\sigma$} if
\[\Delta([a,b])=[\Delta(a),b]+[a,\sigma(b)]\quad (a,b\in \mathcal{A}).\]
A linear map $\xi:\mathcal{A} \rightarrow \mathcal{A}$ is said to be a \textit{Lie triple derivation} if
\[\xi([[a,b],c])=[[\xi(a),b],c]+[[a,\xi(b)],c]+[[a,b],\xi(c)]\quad (a,b\in \mathcal{A}).\]
A linear map $\Lambda:\mathcal{A} \rightarrow \mathcal{A}$ is said to be a \textit{generalized Lie triple derivation associated with the Lie triple derivation $\xi$} if
\[\Lambda([[a,b],c])=[[\Lambda(a),b],c]+[[a,\xi(b)],c]+[[a,b],\xi(c)]\quad (a,b\in \mathcal{A}).\]
Every derivation is a Lie derivation, and every Lie derivation is a generalized Lie derivation. Obviously, Lie derivations are Lie triple derivations, and Lie triple derivations are generalized Lie triple derivations. However, the converse is not true in general. These mappings are very important subjects in the research of Lie structure of algebras. Extensive studies have been performed to characterize these maps on different algebras, and here, for instance, we refer to \cite{beh, ben0, ben, ben2, che, mie, qi1, xia} and the references therein. By routine verifications, it can be seen that $\Delta$ is a generalized Lie derivation associated with the Lie derivation $\sigma$ if and only if $\Delta-\sigma$ is a Lie centralizer. In the following note, we examine the relationship between Lie triple derivations and generalized Lie triple derivations.
\begin{rem}\label{RE}
Let $ \mathcal{A} $ be an algebra. The linear map $\Lambda:\mathcal{A} \rightarrow \mathcal{A}$ is a generalized Lie triple derivation associated with the Lie triple derivation $\xi :\mathcal{A} \rightarrow \mathcal{A}$ if and only if $\Lambda-\xi$ satifies
\[(\Lambda-\xi)([[a,b],c])=[[(\Lambda-\xi)(a),b],c] \]
for all $a,b,c\in \mathcal{A}$. The reason for this is as follows. Suppose the previous identity is established. Set $\phi =\Lambda -\xi$. So
\begin{equation*}
\begin{split}
\Lambda([[a,b],c])&=\xi([[a,b],c])+\phi([[a,b],c]) \\ &
=[[\xi(a),b],c]+[[a,\xi(b)],c]+[[a,b],\xi(c)]+[[\phi(a),b],c] \\&
=[[\Lambda(a),b],c]+[[a,\xi(b)],c]+[[a,b],\xi(c)]
\end{split}
\end{equation*}
for all $a,b,c\in \mathcal{A}$. Hence $\Lambda$ is a generalized Lie triple derivation associated with the Lie triple derivation $\xi$. The converse is clear.
\end{rem}
With regard to the above and Remark \ref{RE}, the idea of defining the Lie triple centralizer is obtained. Let $\phi : \mathcal{A}\rightarrow \mathcal{A} $ be a linear map. We say that $\phi$ is a \textit{Lie triple centralizer} if 
\[ \phi([[a,b],c])=[[\phi(a),b],c] \]
for all $a,b,c\in \mathcal{A}$. It is easily checked that $\phi$ is a Lie triple centralizer on $\mathcal{A}$ if and only if $\phi([[a,b],c])=[[a, \phi(b)],c] $ for all $a,b,c\in \mathcal{A}$. It is clear each Lie centralizer is a Lie triple centralizer, but the converse is not true in general (see Example \ref{EX}). Therefore, the concept of Lie triple centralizer generalizes the concept of Lie centralizer. For $a,b\in \mathcal{A}$, the Jordan product of $a$ and $b$ is denoted by $a\circ b$ and defined by $a\circ b = ab+ba$. A linear mapping $ \psi : \mathcal{A} \to \mathcal{A} $ is called a \textit{Jordan centralizer} if
\[ \psi(a\circ b)=\psi(a)\circ b\] 
for all $a,b\in \mathcal{A}$. The formula $[[a,b], c] = a \circ (b \circ c)- b \circ (a \circ c)$ implies that every Jordan centralizer is also a Lie triple centralizer. By Remark \ref{RE}, $\Lambda$ is a generalized Lie triple derivation associated with the Lie triple derivation $\xi$ if and only if $\Lambda-\xi$ is a Lie triple centralizer. So on an algebra, if we characterize Lie triple centralizers and Lie triple derivations, then we get the characterization of generalized Lie triple derivations. It should be noted that based on this idea, generalized Lie derivations of unital algebras with idempotents in \cite{ben2} have been characterized. The above shows the significance of the definition and characterization of Lie triple centralizers. Suppose that $ \lambda \in Z(\mathcal{A}) $, where $ Z(\mathcal{A}) $ is the center of $ \mathcal{A} $, and $\chi : \mathcal{A} \to Z(\mathcal{A}) $ is a linear mapping which annihilates all second commutators in $ \mathcal{A} $, i.e., $ \chi ( [ [ a,b ] , c ] ) = 0$ for all $ a,b,c\in \mathcal{A} $. 
In this case, the linear mapping $ \phi : \mathcal{A} \to \mathcal{A} $ defined by $ \phi (a) = \lambda a + \chi (a) $ is a Lie triple centralizer, which is called the \textit{proper Lie triple centralizer}. The following example shows that in general every Lie triple centralizer is not necessarily a proper Lie triple centralizer.
\begin{exm} \label{EX}
Let $M_{2}(\mathcal{A}) $ be the algebra of all $2\times 2$ matrices over the algebra 
\[ \mathcal{A} = \left\lbrace \begin{bmatrix}
0 & k_1 & p \\
0 & 0 & k_2 \\
0 & 0 & 0 
\end{bmatrix} : k_1 , k_2 , p \in \mathbb{C} \right\rbrace . \]
Then, the mapping $ \phi :M_{2}(\mathcal{A}) \to M_{2}(\mathcal{A})$ defined by 
\[ \begin{bmatrix}
0 & k_1 & p & 0 & a_1 & b \\
0 & 0 & k_2 & 0 & 0 & a_2 \\
0 & 0 & 0 & 0 & 0 & 0 \\
0 & c_1 & d & 0 & r_1 & q \\
0 & 0 & c_2 & 0 & 0 & r_2 \\
0 & 0 & 0 & 0 & 0 & 0 
\end{bmatrix} \longmapsto \begin{bmatrix}
0 & r_1 & q & 0 & 0 & 0 \\
0 & 0 & r_2 & 0 & 0 & 0 \\
0 & 0 & 0 & 0 & 0 & 0 \\
0 & 0 & 0 & 0 & k_1 & p \\
0 & 0 & 0 & 0 & 0 & k_2 \\
0 & 0 & 0 & 0 & 0 & 0 
\end{bmatrix} , \]
is a Lie triple centralizer which is not a Lie centralizer. In fact, for any $ A, B , C \in M_{2}(\mathcal{A}) $, we have
\[ \phi ( [ [ A, B ] , C ] ) = 0 = [[ \phi (A) , B] , C ] = [ [ A, \phi (B) ] , C ] . \]
If we choose 
\[ A_0= \begin{bmatrix}
0 & 1 & 0 & 0 & 0 & 0 \\
0 & 0 & 1 & 0 & 0 & 0 \\
0 & 0 & 0 & 0 & 0 & 0 \\
0 & 0 & 0 & 0 & 2 & 0 \\
0 & 0 & 0 & 0 & 0 & 1 \\
0 & 0 & 0 & 0 & 0 & 0 
\end{bmatrix} ~ ~ ~ ~ \text{and} ~~~~ B_0= \begin{bmatrix}
0 & 1 & 0 & 0 & 0 & 0 \\
0 & 0 & 1 & 0 & 0 & 0 \\
0 & 0 & 0 & 0 & 0 & 0 \\
0 & 0 & 0 & 0 & 1 & 0 \\
0 & 0 & 0 & 0 & 0 & 2 \\
0 & 0 & 0 & 0 & 0 & 0 
\end{bmatrix}, \]
then it can be easily seen thet $ \phi ( [ A_0, B_0 ] ) \neq [ \phi (A_0) , B_0] $, that is, $ \phi $ is not a Lie centralizer. Also, $\phi$ is not a proper Lie triple centralizer. Because if we assume $\phi$ is a proper Lie triple centralizer, then there exist a $ \lambda \in Z(M_{2}(\mathcal{A}) ) $ and a linear mapping $\chi : M_{2}(\mathcal{A}) \to Z(M_{2}(\mathcal{A})) $ such that $ \phi (A) = \lambda A + \chi (A) $ for all $A \in M_{2}(\mathcal{A})$. By routine verifications, it can be seen that $ Z(M_{2}(\mathcal{A}) )= M_{2}(\mathcal{C})$ where 
\[ \mathcal{C} = \left\lbrace \begin{bmatrix}
0 & 0 & p \\
0 & 0 & 0 \\
0 & 0 & 0 
\end{bmatrix} : p \in \mathbb{C} \right\rbrace . \]
So for any $A\in M_{2}(\mathcal{A})$ and $\gamma\in Z(M_{2}(\mathcal{A}) )$ we have $\gamma A =0$. Hence 
\[\chi (A_0) =\phi (A_0)-\lambda A_0=\phi (A_0)= \begin{bmatrix}
0 & 2 & 0 & 0 & 0 & 0 \\
0 & 0 & 1 & 0 & 0 & 0 \\
0 & 0 & 0 & 0 & 0 & 0 \\
0 & 0 & 0 & 0 & 1 & 0 \\
0 & 0 & 0 & 0 & 0 & 1 \\
0 & 0 & 0 & 0 & 0 & 0 
\end{bmatrix}. \] 
$\chi (A_0)$ is not in $Z(M_{2}(\mathcal{A}) )$ and this is a contradiction. 
\end{exm}
Now the question arises whether any Lie triple centralizer on an algebra is a proper Lie triple centralizer? In this article, we will address this question on generalized matrix algebras and their applications.
\par 
Let $\mathcal{U}$ be a unital generalized matrix algebra (will be defined in the next section), and $ \phi : \mathcal{U} \to \mathcal{U} $ be a Lie triple centralizer. In this paper, we first determine the structure of $ \phi $ (Theorem \ref{t1}). Then, under the appropriate annihilating conditions on $ \mathcal{U} $, we obtain the necessary and sufficient conditions for $ \phi $ to be proper (Theorem \ref{t2}) and from Theorem \ref{t2}, we get sufficient conditions for a Lie triple centralizer to be proper (Corollary \ref{c36}). We then apply our results to characterizations of generalizerd Lie triple derivations on unital generalized matrix algebras (Theorem \ref{tg}). Also, we refer to applications of our results to triangular algebras (Remark \ref{tri}) and some other algebras (Remark \ref{oa}). In Section 2, we set out the preliminaries and tools required. Section 3 is devoted to characterizations of Lie triple centralizers on generalized matrix algebras. In Section 4, we refer to some applications of the results obtained.
\section{Preliminaries and Tools}
A Morita context consists of two algebras $ \mathcal{A} $ and $ \mathcal{B} $, two bimodules $ _{\mathcal{A}}\mathcal{M}_{\mathcal{B}} $ and $ _{\mathcal{B}}\mathcal{N}_{\mathcal{A}} $, and two bimodule homomorphisms called the pairings $ \zeta_{\mathcal{M} \mathcal{N}} : \mathcal{M} \underset{\mathcal{B}}\otimes \mathcal{N} \to \mathcal{A} $ and $ \psi_{\mathcal{N} \mathcal{M}} : \mathcal{N} \underset{\mathcal{A}}\otimes \mathcal{M} \to \mathcal{B} $ satisfying the following commutative diagrams:
\begin{displaymath}
\xymatrix{ \mathcal{M} \underset{\mathcal{B}}\otimes \mathcal{N} \underset{\mathcal{A}}\otimes \mathcal{M} \ar[rr]^{\zeta_{\mathcal{M} \mathcal{N}} \otimes \mathcal{I}_\mathcal{M}} \ar[dd]^{\mathcal{I}_\mathcal{M} \otimes \psi_{\mathcal{N} \mathcal{M}}} & & \mathcal{A} \underset{\mathcal{A}}\otimes \mathcal{M}\ar[dd]^{\cong} \\
& & \\
\mathcal{M} \underset{\mathcal{B}}\otimes \mathcal{B} \ar[rr]^{\cong} & & \mathcal{M} } 
\quad \text{and} \quad
\xymatrix{ \mathcal{N} \underset{\mathcal{A}}\otimes \mathcal{M} \underset{\mathcal{B}}\otimes \mathcal{N} \ar[rr]^{\psi_{\mathcal{N} \mathcal{M}} \otimes \mathcal{I}_\mathcal{N}} \ar[dd]^{\mathcal{I}_\mathcal{N} \otimes \zeta_{\mathcal{M} \mathcal{N}}} & & \mathcal{B} \underset{\mathcal{B}}\otimes \mathcal{N}\ar[dd]^{\cong} \\
& & \\
\mathcal{N} \underset{\mathcal{A}}\otimes \mathcal{A} \ar[rr]^{\cong} & & \mathcal{N} } 
\end{displaymath} 
If $ ( \mathcal{A} , \mathcal{B} , \mathcal{M} , \mathcal{N} , \zeta_{\mathcal{M} \mathcal{N}} , \psi_{\mathcal{N} \mathcal{M}} ) $ is a Morita context, then, the set
\[ \mathcal{U} = \begin{bmatrix}
\mathcal{A} & \mathcal{M} \\
\mathcal{N} & \mathcal{B}
\end{bmatrix} = \left\lbrace \begin{bmatrix}
a & m \\
n & b
\end{bmatrix} : a \in \mathcal{A}, m \in \mathcal{M} , n \in \mathcal{N} , b \in \mathcal{B} \right\rbrace \]
is an algebra under the usual matrix operations. Such an algebra is called a generalized matrix algebra. This type of algebra 
was first introduced by Morita \cite{mor}. If $ \mathcal{A} $ is a unital associative algebra with the identity $1$ that has the 
idempotent $e$ ($e^2=e$), which is non-trivial ($ e \neq 0 ,1 $), then Pierce decomposes of $ \mathcal{A} $ is $ 
\mathcal{A} =e \mathcal{A} e \oplus e \mathcal{A} f \oplus f \mathcal{A} e \oplus f \mathcal{A} f $, where $ f = 1 - e $.
In fact, $ \mathcal{A} $ is a generalized matrix algebra in the form of $ \mathcal{A} = \begin{bmatrix}
e \mathcal{A} e & e \mathcal{A} f \\
f \mathcal{A} e & f \mathcal{A} f
\end{bmatrix} $. So any associative algebra containing a non-trivial idempotent is a generalized algebra. With this in mind, we see that generalized matrix algebra includes a wide range of algebras. In particular, algebra $ M_n(\mathcal{A}) $ of matrices 
$ n \times n $ is a generalized matrix algebra over the unital algebra $ \mathcal{A} $. If $ \mathcal{N} =0 $, then the 
generalized matrix algebra is called the triangular algebra, which is $ \begin{bmatrix}
\mathcal{A} & \mathcal{M} \\
0 & \mathcal{B}
\end{bmatrix} $, and we denote it by $ \mathcal{T} = Tri (\mathcal{A} , \mathcal{M} , \mathcal{B} ) $. Triangular algebras
also include a wide range of algebras, such as nest algebras on Hilbert spaces and upper triangular matrix algebras.
\par 
Suppose that $ \mathcal{A} $ and $ \mathcal{B} $ are unital algebras with identities $1_\mathcal{A}$ and 
$1_\mathcal{B}$, respectively. An $(\mathcal{A}, \mathcal{B})$-bimodule $ \mathcal{M} $ is called unital if 
$1_\mathcal{A} m = m = m 1_\mathcal{B}$, for any $ m \in \mathcal{M} $. It proved that the generalized matrix 
algebra $ \mathcal{U} $ is unital if and only if $ \mathcal{A} $ and $ \mathcal{B} $ are unital algebras with identities $1_\mathcal{A}$ and $1_\mathcal{B}$, respectively, and $ \mathcal{M} $ is a unital $(\mathcal{A}, \mathcal{B})$-bimodule and $ \mathcal{N} $ is a unital $(\mathcal{B}, \mathcal{A})$-bimodule. Then $ I= \begin{bmatrix}
1_\mathcal{A} & 0 \\
0 & 1_\mathcal{B}
\end{bmatrix} $ is the identity element of $ \mathcal{U} $. Given this it is clear that the triangular algebra $ \mathcal{T} = Tri (\mathcal{A} , \mathcal{M} , \mathcal{B} ) $ is unital if and only if $ \mathcal{A} $ and $ \mathcal{B} $ are unital algebras and $ \mathcal{M} $ is a unital $(\mathcal{A}, \mathcal{B})$-bimodule. $ \mathcal{M} $ is called faithful $(\mathcal{A}, \mathcal{B})$-bimodule if $a \mathcal{M} = 0$ ($a \in \mathcal{A}$), then $a = 0$ and from $\mathcal{M} b = 0$ ($b \in \mathcal{B}$) we can deduce that $b = 0$. In studies on generalized matrix algebra, $ \mathcal{M} $ or $ \mathcal{N} $ are often considered as faithful bimodules. But in this article, we propound the following weaker condition on $ \mathcal{U} $:
\[ a \in \mathcal{A}, a \mathcal{M} =0 ~~~ \text{and} ~~~ \mathcal{N} a = 0 \Longrightarrow a = 0 , \]
\[ b \in \mathcal{B} , \mathcal{M} b =0 ~~~ \text{and} ~~~ b \mathcal{N} = 0 \Longrightarrow b = 0 . \]
This weaker condition seems to have been considered at the outset in \cite{ben, gh}. Clearly, if $ \mathcal{M} $ or $ \mathcal{N} $ are faithful bimodules, then the above conditions is established. If $ \mathcal{T} = Tri ( \mathcal{A} , \mathcal{M} , \mathcal{B} ) $ is a triangular algebra, then the annihilating conditions assumed above are equivalent to the faithfulness of $ \mathcal{M} $ as the $(\mathcal{A}, \mathcal{B})$-bimodule.
\par 
Next $ \mathcal{U}$ always indicate a generalized matrix algebra $ \begin{bmatrix}
\mathcal{A} & \mathcal{M} \\
\mathcal{N} & \mathcal{B}
\end{bmatrix} $. Define projection mappings $ \pi_\mathcal{A} : \mathcal{U} \to \mathcal{A} $ and $ \pi_\mathcal{B} : \mathcal{U} \to \mathcal{B} $ by
\[ \pi_\mathcal{A} \left( \begin{bmatrix}
a & m \\
n & b 
\end{bmatrix} \right) = a \quad \text{and} \quad \pi_\mathcal{B} \left( \begin{bmatrix}
a & m \\
n & b 
\end{bmatrix} \right) = b . \]
The following result is proved in \cite{ben}.
\begin{lem}\label{l24} $($\cite[Proposition 2.1]{ben}$)$
Soppose that the unital generalized matrix algebra $ \mathcal{U} $ satisfies
\[ a \in \mathcal{A}, a \mathcal{M} =0 ~~~ \text{and} ~~~ \mathcal{N} a = 0 \Longrightarrow a = 0 , \]
\[ b \in \mathcal{B} , \mathcal{M} b =0 ~~~ \text{and} ~~~ b \mathcal{N} = 0 \Longrightarrow b = 0 . \]
The center of $\mathcal{U}$ is 
\[ Z ( \mathcal{U} ) = \left\lbrace \begin{bmatrix}
a & 0 \\
0 & b
\end{bmatrix} : am = mb , na = bn, \forall m \in \mathcal{M} , \forall n \in 
\mathcal{N} \right\rbrace . \]
Furthermore, there exists a unique algebra isomorphism $ \eta : \pi_\mathcal{A} ( Z (\mathcal{U})) \to \pi_\mathcal{B} ( Z ( 
\mathcal{U} )) $ such that $ a m = m \eta(a) $ and $ n a = \eta (a) n $ for all $ m \in \mathcal{M} $, $ n \in \mathcal{N}$.
\end{lem}
According to Lemma \ref{l24}, we have $ \pi_{\mathcal{A}} ( Z ( \mathcal{U} )) \subseteq Z ( \mathcal{A} ) $ and $ 
\pi_{\mathcal{B}} ( Z ( \mathcal{U} )) \subseteq Z ( \mathcal{B} ) $. In fact, $\pi_{\mathcal{A}} ( Z ( \mathcal{U} ))$ is a 
subalgebra of $ Z ( \mathcal{A} ) $ and $ \pi_{\mathcal{B}} ( Z ( \mathcal{U} )) $ is a subjunctive of $ Z ( \mathcal{B} ) $. In the event that $ \mathcal{T} = Tri ( \mathcal{A}, \mathcal{M} , \mathcal{B} ) $ is a unital triangular algebra that $ \mathcal{M} $ is a faithful ($ \mathcal{A}, \mathcal{B} ) $-bimodulus, then the statement \ref{l24} holds for $ \mathcal{T} $ as well.
\par 
Let $ \mathcal{C} $ be a subset of an algebra $ \mathcal{A} $, we denote by $ \mathcal{C}' $ the commutant of $ \mathcal{C} $ in $ \mathcal{A} $, where
\[ \mathcal{C}' = \{ a \in \mathcal{A} : ac = ca ~\text{for~ any}~ c \in \mathcal{C} \} . \]
\section{characterizations of Lie triplel centralizers}
In this section, we present the main results of this paper. Throughout this section, we assume that $ \mathcal{U} = \begin{bmatrix}
\mathcal{A} & \mathcal{M} \\
\mathcal{N} & \mathcal{B}
\end{bmatrix} $ is a unital generalized matrix algebra. In the following theorem, we first give the structure of Lie triple centralizers on a unital generalized matrix algebras. 
\begin{thm}\label{t1}
A linear map $ \phi : \mathcal{U} \to \mathcal{U} $ is a Lie triplel centralizer if and only if $ \phi $ has the form
\[ \phi \left( \begin{bmatrix}
a & m \\
n & b
\end{bmatrix} \right) = \begin{bmatrix}
\alpha_1 ( a) + \beta_1 (b) & \tau_2 (m) \\
\gamma_3 (n ) & \alpha_4 (a) + \beta_4 (b) 
\end{bmatrix} \]
where $ \alpha_1 : \mathcal{A} \to \mathcal{A} $, $ \beta_1 : \mathcal{B} \to [ \mathcal{A} , \mathcal{A} ]' $, $ \tau_2 : \mathcal{M} \to \mathcal{M} $, $ \alpha_4 : \mathcal{A} \to [ \mathcal{B} , \mathcal{B} ]' $, $ \beta_4 : \mathcal{B} \to \mathcal{B} $ and $ \gamma_3 : \mathcal{N} \to \mathcal{N} $ are linear mappings satisfying the following conditions:
\begin{enumerate}
\item[(i)]
$ \alpha_1 $ is a Lie triplel centralizer on $ \mathcal{A} $, $ [ \alpha_4 (a) , b_1 ], b_2 ] = 0 $, $ \alpha_4 ([ [a_1 , a_2 ] , a_3 ] ) = 0 $ and 
\[ \alpha_1 ( mn ) - \beta_1 (nm) = \tau_2(m) n = m \gamma_3 (n), \quad (m \in \mathcal{M} , n\in \mathcal{N}). \]
\item[(ii)]
$ \beta_4 $ is a Lie triplel centralizer on $ \mathcal{B} $, $ [ [ \beta_1 (b) , a_1 ] , a_2 ] = 0 $, $ \beta_1 ( [ [ b_1 , b_2 ] , b_3 ] ) = 0 $ and 
\[ \alpha_4 ( nm ) - \beta_4 (mn) = n \tau_2(m) = \gamma_3 (n) m , \quad (m \in \mathcal{M} , n\in \mathcal{N}). \]
\item[(iii)]
$ \tau_2 (am) = a \tau_2 (m) = \alpha_1 (a) m - m \alpha_4 (a) $ and $ \tau_2 (mb) = \tau_2 (m) b = m \beta_4 (b) - \beta_1(b) m $ for any $ a \in \mathcal{A} $, $ b \in \mathcal{B} $ and $ m \in \mathcal{M} $.
\item[(iv)]
$ \gamma_3 (na) = \gamma_3 (n) a = n \alpha_1 (a) - \alpha_4 (a) n $ and $ \gamma_3 (b n) = b \gamma_3 (m) = \beta_4 (b) n - n \beta_1(b) $ for any $ a \in \mathcal{A} $, $ b \in \mathcal{B} $ and $ m \in \mathcal{M} $.
\end{enumerate}
\end{thm}
\begin{proof}
Suppose that the map $ \phi $ on $ \mathcal{U} $ has the form
\[ \phi \left( \begin{bmatrix}
a & m \\
n & b
\end{bmatrix} \right) = \begin{bmatrix}
\alpha_1(a) + \beta_1 (b) + \tau_1(m) + \gamma_1(n) & \alpha_2 (a) + \beta_2 (b) + \tau_2(m) + \gamma_2(n) \\
\alpha_3(a) + \beta_3 (b) + \tau_3(m) + \gamma_3(n) & \alpha_4(a) + \beta_4 (b) + \tau_4 (m) + \gamma_4(n)
\end{bmatrix} \]
where 
\[ \alpha_1 : \mathcal{A} \to \mathcal{A}, ~~ \alpha_2 : \mathcal{A} \to \mathcal{M}, ~~ \alpha_3 : \mathcal{A} \to \mathcal{N}, ~~ \alpha_4 : \mathcal{A} \to \mathcal{B}, \]
\[ \beta_1 : \mathcal{B} \to \mathcal{A}, ~~ \beta_2 : \mathcal{B} \to \mathcal{M}, ~~ \beta_3 : \mathcal{B} \to \mathcal{N}, ~~ \beta_4 : \mathcal{B} \to \mathcal{B}, \]
\[ \tau_1 : \mathcal{M} \to \mathcal{A}, ~~ \tau_2 : \mathcal{M} \to \mathcal{M}, ~~ \tau_3 : \mathcal{M} \to \mathcal{N}, ~~ \tau_4 : \mathcal{M} \to \mathcal{B}, \]
\[ \gamma_1 : \mathcal{N} \to \mathcal{A}, ~~ \gamma_2 : \mathcal{N} \to \mathcal{M}, ~~ \gamma_3 : \mathcal{N} \to \mathcal{N}, ~~ \gamma_4 : \mathcal{NB} \to \mathcal{B}, \]
are linear mappings.
\par 
We assume that $ \phi $ is a Lie triplel centralizer. So for any $ A, B, C \in \mathcal{U} $, we consider
\begin{equation}\label{eq1}
\phi ( [ [A , B ] , C ] ) = [ [ \phi (A) , B ] , C ] . 
\end{equation}
By taking $ A= \begin{bmatrix}
a & 0 \\
0 & 0
\end{bmatrix} $, $ B= \begin{bmatrix}
0 & 0 \\
0 & b_1 
\end{bmatrix} $ and $ C= \begin{bmatrix}
0 & 0 \\
0 & b_2
\end{bmatrix} $ in \eqref{eq1}, we deduce that 
\begin{align*}
0 = \phi ( [ [A , B ] , C ] ) 
= [ [ \phi (A) , B ] , C ] = 
\begin{bmatrix}
0 & \alpha_2 (a) b_1 b_2 \\
b_2 b_1 \alpha_3 (a) & [[ \alpha_4(a) , b_1 ] , b_2 ] 
\end{bmatrix}
\end{align*}
So 
\[ \alpha_2 (a) b_1 b_2 = 0 , \quad b_2 b_1 \alpha_3 (a)= 0, \quad [[ \alpha_4(a) , b_1 ] , b_2 ] =0 , \]
for any $ a \in \mathcal{A} $ and any $ b_1 , b_2 \in \mathcal{B} $. If we set $ b_1 = b_2 = 1_\mathcal{B} $, then $ \alpha_2 = 0 $ and $ \alpha_3 = 0 $. 
\par 
Let $ A= \begin{bmatrix}
a_1 & 0 \\
0 & 0
\end{bmatrix} $, $ B= \begin{bmatrix}
a_2 & 0 \\
0 & 0 
\end{bmatrix} $ and $ C= \begin{bmatrix}
a_3 & 0 \\
0 & 0
\end{bmatrix} $ in \eqref{eq1}, to obtain 
\begin{align*}
\begin{bmatrix}
\alpha_1 ( [[a_1, a_2] , a_3 ] ) & 0 \\
0 & \alpha_4 ( [[a_1, a_2] , a_3 ] ) 
\end{bmatrix} & = \phi ( [ [A , B ] , C ] ) \\
& = [ [ \phi (A) , B ] , C ] \\
& =
\begin{bmatrix}
( [ \alpha_1 (a_1) , a_2 ] ) a_3 - a_3 ( [ \alpha_1 (a_1) , a_2 ] ) & 0 \\
0 & 0 
\end{bmatrix}
\end{align*}
from the above relation, we see that $ \alpha_1 $ is a Lie triple centralizer on $ \mathcal{A} $ and 
\[\alpha_4 ( [[a_1, a_2] , a_3 ] ) =0. \]
\par 
If 
$ A= \begin{bmatrix}
0 & 0 \\
0 & b
\end{bmatrix} $, $ B= \begin{bmatrix}
a_1 & 0 \\
0 & 0 
\end{bmatrix} $ and $ C= \begin{bmatrix}
a_2 & 0 \\
0 & 0
\end{bmatrix} $, then matrix equation \eqref{eq1} implies
\begin{align*}
0 = \phi ( [ [A , B ] , C ] ) 
= [ [ \phi (A) , B ] , C ] = 
\begin{bmatrix}
[[ \beta_1(b) , a_1 ] , a_2 ] & \beta_2 (b) a_1 a_2 \\
a_2 a_1 \beta_3 (b) & 0
\end{bmatrix}
\end{align*}
Thus
\[ [[ \beta_1(b) , a_1 ] , a_2 ]=0, \quad \beta_2 (b) a_1 a_2 = 0 , \quad a_2 a_1 \beta_3 (b) = 0 , \]
for any $ b \in \mathcal{B} $ and any $ a_1 , a_2 \in \mathcal{A} $. Now putting $ a_1 = a_2 = 1_\mathcal{A} $ in above, we obtain $ \beta_2 = 0 $ and $ \beta_3 = 0 $. Similarly, by taking $ A= \begin{bmatrix}
0 & 0 \\
0 & b_1
\end{bmatrix} $, $ B= \begin{bmatrix}
0 & 0 \\
0 & b_2
\end{bmatrix} $ and $ C= \begin{bmatrix}
0 & 0 \\
0 & b_3
\end{bmatrix} $, we can find that $ \beta_4 $ is a Lie triple centralizer on $ \mathcal{B} $ and $ \beta_1 ( [[b_1, b_2] , b_3 ] ) =0 $. 

Let us choose $ A= \begin{bmatrix}
a & 0 \\
0 & 0
\end{bmatrix} $, $ B= \begin{bmatrix}
0 & m \\
0 & 0 
\end{bmatrix} $ and $ C= \begin{bmatrix}
1_\mathcal{A} & 0 \\
0 & 0
\end{bmatrix} $ in \eqref{eq1}. Then, we find that 
\begin{align*}
\begin{bmatrix}
-\tau_1 (a m) & - \tau_2 (am) \\
- \tau_3 (am) & \tau_4 (am ) 
\end{bmatrix} & = \phi ( [ [A , B ] , C ] ) \\
& = [ [ \phi (A) , B ] , C ] = 
\begin{bmatrix}
0 & -( \alpha_1 (a) m - m \alpha_4 (am) \\
0 & 0 
\end{bmatrix}
\end{align*}
Thus 
\[ \tau_1 (am ) = 0, \quad \tau_3 (am ) = 0 , \quad \tau_4 (am ) = 0 , \]
for any $ a \in \mathcal{A} $ and any $ m \in \mathcal{M} $. Taking $ a= 1_\mathcal{A} $ in above, we see that $ \tau_1 = 0 $, $ \tau_3 = 0 $ and $ \tau_4 = 0 $. Also 
\begin{equation}\label{eq2}
\tau_2 (am) = \alpha_1 (a) m - m \alpha_4 (a)
\end{equation}
for any $ a \in \mathcal{A} $ and any $ m \in \mathcal{M} $. Similarly, by taking $ A= \begin{bmatrix}
0 & m \\
0 & 0
\end{bmatrix} $, $ B= \begin{bmatrix}
a & 0 \\
0 & 0 
\end{bmatrix} $ and $ C= \begin{bmatrix}
0 & 0 \\
0 & 1_\mathcal{B}
\end{bmatrix} $, we have $ \tau_2 (am ) = a \tau_2 (m) $ for any $ a \in \mathcal{A} $ and any $ m \in \mathcal{M} $. Furthermore, by taking $ A= \begin{bmatrix}
0 & m \\
0 & 0
\end{bmatrix} $, $ B= \begin{bmatrix}
0 & 0 \\
0 & b
\end{bmatrix} $ and $ C= \begin{bmatrix}
0 & 0 \\
0 & 1_\mathcal{B}
\end{bmatrix} $ and also 
$ A= \begin{bmatrix}
0 & 0 \\
0 & b
\end{bmatrix} $, $ B= \begin{bmatrix}
0 & 0 \\
0 & 1_\mathcal{B}
\end{bmatrix} $ and $ C= \begin{bmatrix}
0 & 0 \\
0 & 1_\mathcal{B}
\end{bmatrix} $, we have
\[ \tau_2 ( m b) = \tau_2 (m) b , \quad \quad \tau_2 (mb) = m \beta_4 (b) - \beta_1 ( b) m , \]
for any $ b \in \mathcal{B} $ and any $ m \in \mathcal{M} $.

If we take $ A= \begin{bmatrix}
0 & m \\
0 & 0
\end{bmatrix} $, $ B= \begin{bmatrix}
a & 0 \\
0 & 0 
\end{bmatrix} $ and $ C= \begin{bmatrix}
0 & 0 \\
0 & 1_\mathcal{B}
\end{bmatrix} $ in \eqref{eq1}, then we have 
\begin{equation}\label{eq3}
\phi ( [ [A , B ] , C ] ) = \begin{bmatrix}
\gamma_1 (n a) & \gamma_2 ( na) \\
\gamma_3 (na) & \gamma_4 (na) 
\end{bmatrix}
\end{equation}
and 
\begin{equation}\label{eq4}
[ [ \phi (A) , B ] , C ] = \begin{bmatrix}
0 & 0 \\
n \alpha_1 (a) - \alpha_4 (a) n & 0 
\end{bmatrix}
\end{equation}
Combining \eqref{eq3} and \eqref{eq4}, we get 
\[ \gamma_1 (na) = 0 , \quad \gamma_2 (na) = 0 , \quad \gamma_4 (na) = 0 \]
for any $ a \in \mathcal{A} $ and any $ n \in \mathcal{N} $. Set $ a = 1_\mathcal{A} $, we find $ \gamma_1 = 0 $, $ \gamma_2 = 0 $ and $ \gamma_4 = 0 $. Also matrix equations obove, implies 
\[ \gamma_3 ( na) = \alpha_1 (a) n - n \alpha_4 (a) n , \]
for any $ a \in \mathcal{A} $ and any $ n \in \mathcal{N} $. Similarly, by considering $ A= \begin{bmatrix}
0 & 0 \\
n & 0
\end{bmatrix} $, $ B= \begin{bmatrix}
a & 0 \\
0 & 0 
\end{bmatrix} $ and $ C= \begin{bmatrix}
1_\mathcal{A} & 0 \\
0 & 0
\end{bmatrix} $, we obtain $ \gamma_3 (na ) = \gamma_3 ( na) $ for any $ a \in \mathcal{A} $ and any $ n \in \mathcal{N} $. Similarly, by taking $ A= \begin{bmatrix}
0 & 0 \\
n & 0
\end{bmatrix} $, $ B= \begin{bmatrix}
0 & 0 \\
0 & b
\end{bmatrix} $ and $ C= \begin{bmatrix}
1_\mathcal{A} & 0 \\
0 & 0
\end{bmatrix} $ and also 
$ A= \begin{bmatrix}
0 & 0 \\
0 & b
\end{bmatrix} $, $ B= \begin{bmatrix}
0 & 0 \\
n & 0
\end{bmatrix} $ and $ C= \begin{bmatrix}
1_\mathcal{A} & 0 \\
0 & 0
\end{bmatrix} $ in \eqref{eq1}, we arrive at 
\[ \gamma_3 ( b n) = b \gamma_3 (n) , \quad \quad \gamma_3 ( b n ) = \beta_4 (b) n - n \beta_1 (b) , \]
for any $ b \in \mathcal{B} $ and any $ n \in \mathcal{N} $.

Now, by using \eqref{eq1} for $ A= \begin{bmatrix}
1_\mathcal{A} & 0 \\
0 & 0
\end{bmatrix} $, $ B= \begin{bmatrix}
0 & m \\
0 & 0
\end{bmatrix} $ and $ C= \begin{bmatrix}
0 & 0 \\
n & 0
\end{bmatrix} $, we get
\[ \alpha_1 ( mn ) - \beta_1 (nm) = \tau_2 (m) n , \]
\[ \alpha_4 (nm) - \beta_4 ( nm) = n \tau_2 (m) , \]
for any $ m \in \mathcal{M} $ and any $ n \in \mathcal{N} $. Similarly, by putting 
$ A= \begin{bmatrix}
1_\mathcal{A} & 0 \\
0 & 0
\end{bmatrix} $, $ B= \begin{bmatrix}
0 & m \\
0 & 0
\end{bmatrix} $ and $ C= \begin{bmatrix}
0 & 0 \\
n & 0
\end{bmatrix} $ in \eqref{eq1}, we have
\[ \alpha_1 ( mn ) - \beta_1 (nm) = m \gamma_3 (n) , \]
\[ \alpha_4 (nm) - \beta_4 ( nm) = \gamma_3 (n) m , \]
for any $ m \in \mathcal{M} $ and any $ n \in \mathcal{N} $.

Conversely, suppose that $ \phi $ is a linear map on $ \mathcal{U} = \begin{bmatrix}
\mathcal{A} & \mathcal{M} \\
\mathcal{N} & \mathcal{B}
\end{bmatrix} $ has the form 
\[ \phi \left( \begin{bmatrix}
a & m \\
n & b
\end{bmatrix} \right) = \begin{bmatrix}
\alpha_1 ( a) + \beta_1 (b) & \tau_2 (m) \\
\gamma_3 (n ) & \alpha_4 (a) + \beta_4 (b) 
\end{bmatrix} \] 
and satisfies the assumptions (i), (ii), (iii) and (iv). Then, it is easy to check that $ \phi $ satisfies the relation 
\[ \phi ( [ [ A,B] , C ] = [ [ \phi (A), B ] , C ] , \quad \quad A, B, C \in \mathcal{U} . \]
This completes the proof.
\end{proof}
In the next Corollary, we consider the annihilating conditions on $\mathcal{U}$ and in this case we see that the conditions $ \alpha_4 ( [ [ a_1 , a_2 ] , a_3 ] ) = 0 $ and $ \beta_1 ( [ [ b_1 , b_2 ] , b_3 ] ) = 0 $ in the above theorem can be omitted, and also $ \alpha_4 (a) \in Z( \mathcal{B} ) $ 
and $ \beta_1 (b) \in Z( \mathcal{A} ) $.
\begin{cor}\label{c2}
Let $ \mathcal{U} $ satisfies 
\[ a \in \mathcal{A}, a \mathcal{M} =0 ~~~ \text{and} ~~~ \mathcal{N} a = 0 \Longrightarrow a = 0 , \]
\[ b \in \mathcal{B} , \mathcal{M} b =0 ~~~ \text{and} ~~~ b \mathcal{N} = 0 \Longrightarrow b = 0 . \]
A linear map $ \phi : \mathcal{U} \to \mathcal{U} $ is a Lie triplel centralizer if and only if $ \phi $ has the form
\[ \phi \left( \begin{bmatrix}
a & m \\
n & b
\end{bmatrix} \right) = \begin{bmatrix}
\alpha_1 ( a) + \beta_1 (b) & \tau_2 (m) \\
\gamma_3 (n ) & \alpha_4 (a) + \beta_4 (b) 
\end{bmatrix} \]
where $ \alpha_1 : \mathcal{A} \to \mathcal{A} $, $ \beta_1 : \mathcal{B} \to Z( \mathcal{A} ) $, $ \tau_2 : \mathcal{M} \to \mathcal{M} $, $ \alpha_4 : \mathcal{A} \to Z( \mathcal{B} ) $, $ \beta_4 : \mathcal{B} \to \mathcal{B} $ and $ \gamma_3 : \mathcal{N} \to \mathcal{N} $ are linear mappings satisfying the following conditions:
\begin{enumerate}
\item[(i)]
$ \alpha_1 $ is a Lie triplel centralizer on $ \mathcal{A} $, and 
\[ \alpha_1 ( mn ) - \beta_1 (nm) = \tau_2(m) n = m \gamma_3 (n), \quad (m \in \mathcal{M} , n\in \mathcal{N}). \]
\item[(ii)]
$ \beta_4 $ is a Lie triplel centralizer on $ \mathcal{B} $, and 
\[ \alpha_4 ( nm ) - \beta_4 (mn) = n \tau_2(m) = \gamma_3 (n) m , \quad (m \in \mathcal{M} , n\in \mathcal{N}). \]
\item[(iii)]
$ \tau_2 (am) = a \tau_2 (m) = \alpha_1 (a) m - m \alpha_4 (a) $ and $ \tau_2 (mb) = \tau_2 (m) b = m \beta_4 (b) - \beta_1(b) m $ for any $ a \in \mathcal{A} $, $ b \in \mathcal{B} $ and $ m \in \mathcal{M} $.
\item[(iv)]
$ \gamma_3 (na) = \gamma_3 (n) a = n \alpha_1 (a) - \alpha_4 (a) n $ and $ \gamma_3 (b n) = b \gamma_3 (m) = \beta_4 (b) n - n \beta_1(b) $ for any $ a \in \mathcal{A} $, $ b \in \mathcal{B} $ and $ m \in \mathcal{M} $.
\end{enumerate}
\end{cor}
\begin{proof}
Let $ \phi $ be a Lie triplel centralizer, So $ \phi $ satisfies conditions of Theorem \ref{t1}. Hence by Theorem \ref{t1}-(iii), we have
\begin{align*}
\alpha_1 ( & [ [ a_1 , a_2 ] , a_3 ] ) m - m \alpha_4 ( [ [ a_1 , a_2 ] , a_3 ] ) = \tau_2 ( ( [ [ a_1 , a_2 ] , a_3 ] m ) \\
& = \tau_2 \left( ( a_1 a_2 a_3 - a_2 a_1 a_3 - a_3 a_1 a_2 + a_3 a_2 a_1 ) m \right) \\
& = \tau_2 ( (a_1 a_2 a_3) m ) - \tau_2 ( ( a_2 a_1 a_3 ) m ) - \tau_2 ( (a_3 a_1 a_2 ) m ) + \tau_2 ( ( a_3 a_2 a_1 ) m ) \\
& = ( \alpha_1 (a_1) a_2 a_3 ) m - m ( a_2 a_3 \alpha_4 (a_1) ) - ( a_2 \alpha_1 (a_1 ) a_3 ) m + m ( a_2 a_3 \alpha_4 
(a_1) ) \\
& ~~~~ - ( a_3 \alpha_1 (a_1) a_2 ) m + m ( a_3 a_2 \alpha_4 (a_1) ) + ( a_3 a_2 \alpha_1(a_1) ) m - m ( a_3 a_2 
\alpha_4 (a_1) ) \\
& = [ [ \alpha_1 (a_1) , a_2 ] , a_3 ] 
\end{align*}
for any $ a_1 , a_2 , a_3 \in \mathcal{A} $ and $ m \in \mathcal{M} $. By the fact that $ \alpha_1 $ is a Lie triplel centralizer on $ \mathcal{A} $, we have $ \mathcal{M} \alpha_4 ( [ [ a_1 , a_2 ] , a_3 ] ) = 0 $. On the other hand by using Theorem \ref{t1}-(iii), we see that
\begin{align*}
\alpha_1 ( & n [ [ a_1 , a_2 ] , a_3 ] ) - \alpha_4 ( [ [ a_1 , a_2 ] , a_3 ] ) n = \gamma_3 ( n ( [ [ a_1 , a_2 ] , a_3 ]) \\
& = \gamma_3 \left( n ( a_1 a_2 a_3 - a_2 a_1 a_3 - a_3 a_1 a_2 + a_3 a_2 a_1 ) \right) \\
& = \gamma ( n (a_1 a_2 a_3) ) - \gamma_3 ( n ( a_2 a_1 a_3 ) ) - \gamma_3 ( n (a_3 a_1 a_2 ) ) + \gamma_3 ( n ( a_3 a_2 a_1 ) ) \\
& = n ( \alpha_1 (a_1) a_2 a_3 ) - n ( a_2 a_3 \alpha_4 (a_1) ) - n ( a_2 \alpha_1 (a_1 ) a_3 ) + n ( a_2 a_3 \alpha_4 
(a_1) ) \\
& ~~~~ - n( a_3 \alpha_1 (a_1) a_2 ) + n ( a_3 a_2 \alpha_4 (a_1) ) + n ( a_3 a_2 \alpha_1(a_1) ) - n ( a_3 a_2 
\alpha_4 (a_1) ) \\
& = [ [ \alpha_1 (a_1) , a_2 ] , a_3 ] 
\end{align*}
for any $ a_1 , a_2 , a_3 \in \mathcal{A} $ and $ n \in \mathcal{N} $. Since $ \alpha_1 $ is a Lie triplel centralizer on $ \mathcal{A} $, it follows that $ \alpha_4 ( [ [ a_1 , a_2 ] , a_3 ] ) \mathcal{N} = 0 $.
By assumption, we have $ \alpha_4 ( [ [ a_1 , a_2 ] , a_3 ] ) = 0 $ for any $ a_1 , a_2 , a_3 \in \mathcal{A} $. 
Similarly, by using the argument above we have $ \beta_1 ( [ [ b_1 , b_2 ] , b_3 ] ) = 0 $ for any $ b_1 , b_2 , b_3 \in 
\mathcal{B} $. Also for $ a \in \mathcal{A} $, any $ b \in \mathcal{B} $ and any $ m \in \mathcal{M} $, we deduce 
\[ \tau_2 ( a m b ) = \tau_2 ( am ) b = ( \alpha_1 (a) m - m \alpha_4 (a) ) b = \alpha_1 (a) m b - m \alpha_4 (a) b \]
on the other hand
\[ \tau_2 ( a m b) = \alpha_1 (a) m b - m b \alpha_4 (a) . \]
Comparing these identities, we deduce 
\[ m \alpha_4 (a) b = m b \alpha_4 (a) . \]
So 
\[ \mathcal{M} ( \alpha_4 (a) b - b \alpha_4 (a) ) = 0 . \]
Similarly, by using the argument above and condition (iv) we have
\[ ( \alpha_4 (a) b - b \alpha_4 (a) ) \mathcal{N} = 0 . \]
Therefore, $ \alpha_4 (a) \in Z( \mathcal{B} ) $ for $ a \in \mathcal{A} $ as by assumption. By using the argument above for any $ b \in \mathcal{B} $, we see that $ \beta_1 (b) \in Z( \mathcal{A} ) $. The converse is easily checked.
\end{proof}
In the next Theorem, we present the necessary and sufficient conditions for a Lie triple centralizer mapping on $ \mathcal{U} $ to be proper. 
\begin{thm}\label{t2}
Let $ \mathcal{U} $ satisfies the following conditions:
\[ a \in \mathcal{A}, a \mathcal{M} =0 ~~~~ \text{and} ~~~~ \mathcal{N} a = 0 \Rightarrow a = 0 , \]
\[ b\in \mathcal{B}, \mathcal{M} b =0 ~~~~ \text{and} ~~~~ b \mathcal{N} = 0 \Rightarrow b =0 . \]
Suppose that $ \phi : \mathcal{U} \to \mathcal{U} $ is a linear Lie triple centralizer and
\[ \phi \left( \begin{bmatrix}
a & 0 \\
0 & 0
\end{bmatrix} \right) = \begin{bmatrix}
\ast & \ast \\
\ast & \alpha_4 (a)
\end{bmatrix} \]
and
\[ \phi \left( \begin{bmatrix}
0 & 0 \\
0 & b
\end{bmatrix} \right) = \begin{bmatrix}
\beta_1(b) & \ast \\
\ast & \ast
\end{bmatrix} \]
for any $ a \in \mathcal{A} $ and any $ b \in \mathcal{B} $. Then the following statements are equivalent:
\begin{enumerate}
\item[(i)]
$ \phi $ on $ \mathcal{U} $ is a proper Lie triple centralizer; 
\item[(ii)]
$ \alpha_4 (\mathcal{A} ) \subseteq \pi_\mathcal{B} ( Z (\mathcal{U} )) $ and $ \beta_1 (\mathcal{B} ) \subseteq \pi_\mathcal{A} ( Z (\mathcal{U} )) $;
\item[(iii)]
$ \alpha_4 (1_\mathcal{A} ) \in \pi_\mathcal{B} ( Z (\mathcal{U} )) $ and $ \beta_1 (1_\mathcal{B} ) \in \pi_\mathcal{A} ( Z (\mathcal{U} )) $.
\end{enumerate}
\end{thm}
\begin{proof}
According to Theorem \ref{t1}, $ \phi$ has the following form.
\[ \phi \left( \begin{bmatrix}
a & m \\
n & b
\end{bmatrix} \right) = \begin{bmatrix}
\alpha_1 (a) + \beta_1 (b) & \tau_2 (m) \\
\gamma_3 (n) & \alpha_4 (a) + \beta_4 (b) 
\end{bmatrix} \]
where $ \alpha_1$, $ \beta_1 $, $ \tau_2 $, $ \gamma_3 $, $ \alpha_4 $ and $ \beta_4 $ are linear maps with the properties mentioned in Theorem \ref{t1}.\\

$ (i) \Rightarrow (ii) $. Let $ a \in \mathcal{A} $, $ m \in \mathcal{M} $ and $ n \in \mathcal{N} $ are arbitrary elements. We set $ A= \begin{bmatrix}
0 & am \\
na & 0
\end{bmatrix} $. Thus 
\[ \phi (A) = \begin{bmatrix}
0 & \tau_2 (am) \\
\gamma_3 ( na) & 0 
\end{bmatrix} = \begin{bmatrix}
0 & \alpha_1 (a) n - m \alpha_4 (a) \\
n \alpha_1 (a) - \alpha_4 (a) n & 0
\end{bmatrix} . \]
Let $ \lambda = \begin{bmatrix}
a_1 & 0 \\
0 & \eta (a_1)
\end{bmatrix} $ where $ a_1 \in \pi_\mathcal{A} ( Z(\mathcal{U})) $ and $ \chi (A) = \begin{bmatrix}
a_2 & 0 \\
0 & \eta (a_2 )
\end{bmatrix} $ with $ a_2 \in \pi_\mathcal{A} ( Z(\mathcal{U})) $. By assumption, we have
\[ \phi (A) = \lambda A + \chi (A) = \begin{bmatrix}
a_2 & a_1 a m \\
\eta (a_1) n a & \eta(a_2) 
\end{bmatrix} \]
By comparing the identities for $ \phi $, it follows that
\[ \alpha_1 (a) m - m \alpha_4 (a) = a_1 a m , \]
\[ n \alpha_1 (a) - \alpha_4 (a) n = \eta (a_1) n a = n a_1 a . \]
Therefore
\[ (a_1 a - \alpha_1 (a)) m = m \alpha_4 (a) , \]
\[ n ( a_1 a - \alpha_1 (a) ) = \alpha_4 (a) n . \]
Since $ m \in \mathcal{M} $, $ n \in \mathcal{N} $ and $ \alpha_4 (a) \in Z(\mathcal{B} ) $ are arbitrary, hence
$ \alpha_4(a) \in \pi_\mathcal{B} ( Z (\mathcal{U})) $.
Now, for arbitrary elements $ b \in \mathcal{B} $, $ m \in \mathcal{M} $ and $ n \in \mathcal{N} $, we taking
\[ \lambda = \begin{bmatrix}
a_1 & 0 \\
0 & \eta (a_1)
\end{bmatrix} , \quad B = \begin{bmatrix}
0 & m b \\
bn & 0
\end{bmatrix}, \quad \text{and} \quad \chi(B) = \begin{bmatrix}
a_2 & 0 \\
0 & \eta (a_2)
\end{bmatrix} . \]
By the assumption and Theorem \ref{t1} we deduce that 
\[ \phi (B) = \begin{bmatrix}
0 & m \beta_4 (b) - \beta_1 (b) m \\
\beta_4 (b) n - n \beta_1 (b) & 0 
\end{bmatrix} \]
and 
\[ \phi (B) = \begin{bmatrix}
a_2 & a_1 m b \\
\eta(a_1) b n & \eta (a_2 )
\end{bmatrix} \]
So we arrive at 
\[ m ( \beta_4 (b ) - \eta (a_1) b ) = \beta_1 (b) m , \]
\[ ( \beta_4 (b ) - \eta(a_1) b ) n = n \beta_1(b) \]
for any $ m \in \mathcal{M} $ and any $ n \in \mathcal{N} $, and since $ \beta_1 (b) \in Z ( \mathcal{A}) $, we have $ \beta_1(b) \in \pi_\mathcal{A} ( Z( \mathcal{U} )) $. \\

$ (ii) \Rightarrow (i) $. According to Hypothesis, we define the following well-defined functions.
\[ \bar{\alpha} : \mathcal{A} \to \mathcal{A}; \quad \bar{\alpha} (a) = \alpha_1(a) - \eta^{-1} ( \alpha_4(a)) , \]
and 
\[ \bar{\beta} : \mathcal{B} \to \mathcal{B}; \quad \bar{\beta} (b) = \alpha_4(b) - \eta ( \beta_1(b)) , \]
where $ \bar{\alpha} $ and $ \bar{\beta} $ are linear maps and by using Theorem \ref{t1}, we have
\[ \tau_2 (am) = \bar{\alpha} (a) m , ~~ \tau_2(mb) = m \bar{\beta} (a)~~~~ ( a\in \mathcal{A}, b \in \mathcal{B} , m \in \mathcal{M} ) . \]
Thus
\[ \tau_2 (m) = \bar{\alpha} (1_\mathcal{A} ) m = m \bar{\beta} ( 1_\mathcal{B} ) , \]
for any $ m \in \mathcal{M} $. Also, with similar to part $ (iv) $ of Theorem \ref{t1} we arrive at
\[ \gamma_3 (n ) = n \bar{\alpha} (1_\mathcal{A} ) = \bar{\beta} (1_\mathcal{B} ) n \]
for any $ n \in \mathcal{N} $. By Lemma \ref{l24}, it results $ \bar{\alpha} (1_\mathcal{A} ) \in \pi_\mathcal{A} ( Z(\mathcal{U} )) $, $ \bar{\beta } (1_\mathcal{B}) \in \pi_\mathcal{B} ( Z(\mathcal{U} )) $ and $ \eta ( \bar{\alpha} (1_\mathcal{A} ))= \bar{\beta } (1_\mathcal{B} ) $.

For any $ a, a' \in \mathcal{A} $, $ m \in \mathcal{M} $ and $ n \in \mathcal{M} $, we have
\[ a a' \tau_2 (m) = \bar{\alpha} ( a a') m , ~~~ \text{and}~~~ a a' \tau_2(m) = a \alpha(a')m . \]
Therefore
\[ ( \bar{\alpha} (aa') - a \bar{\alpha} (a') ) \mathcal{M} = 0 , \]
also
\[ \gamma_3 ( n a a') = n \bar{\alpha} (a a') , ~~~ \text{and} ~~~ \gamma_3 ( n a a' ) = n a \bar{\alpha} (a' ) . \]
Thus
\[ \mathcal{N} ( \bar{\alpha} (a a') - a \bar{\alpha} (a' ) ) = 0 . \]
Now, by the assumption, we have
\[ \bar{\alpha} (a a') = a \bar{\alpha} (a') \]
for any $ a , a' \in \mathcal{A} $. So 
\[ \bar{\alpha} (a a') = \bar{\alpha} (a 1_\mathcal{A}) = a \bar{\alpha} (1_\mathcal{A}) \]
and since $ \bar{\alpha} (1_\mathcal{A}) \in \pi_\mathcal{A} ( Z(\mathcal{U} )) $, hence
\[ \bar{\alpha} (a ) = a \bar{\alpha} (1_\mathcal{A}) = \bar{\alpha} (1_\mathcal{A}) a , \]
for any $ a \in \mathcal{A} $.

For $ b ,b' \in \mathcal{B} $, $ m \in \mathcal{B} $ and $ n \in \mathcal{N} $, we get
\[ \tau_2 (mb) = \tau_2 (m) b = m \bar{\beta} (b) , \]
\[ \gamma_3 ( b n ) = b \gamma_3 (n) = \bar{\beta} (b) n . \]
According to assumption, identities above, $ \bar{\beta} (1_\mathcal{B} ) \in \pi_\mathcal{B} (Z(\mathcal{U})) $ and the 
same manner before, it results
\[ \bar{\beta} (b) = \bar{\beta} (1_\mathcal{B}) b = b \bar{\beta} (1_\mathcal{B}) \]
for any $ b \in \mathcal{B} $. 

For any $ m \in \mathcal{M} $ and any $ n \in \mathcal{N} $, we have
\begin{align*}
\alpha_1 (mn) - \beta_1 (nm) & = \tau_2 (m)n \\
& = \bar{\alpha} (1_\mathcal{A} ) m n \\
& = \bar{\alpha} ( mn ) \\
&= \alpha_1 (mn) - \eta^{-1} ( \alpha_4 (mn)) 
\end{align*}
So $ \beta_1 ( nm ) = \eta^{-1} ( \alpha_4 (mn)) $ and 
\[ \alpha_1 (mn) - \eta^{-1} ( \alpha_4 (mn)) = \tau_2 (m) n . \]
Thus
\[ \bar{\alpha} ( mn) = \tau_2 (m) n , \]
for any $ m \in \mathcal{M} $ and $ n \in \mathcal{N} $.

According to the obtained result and the equations $ (ii) $ of theorem \ref{t1} and the similar method above, we conclude that
\[ \bar{\beta} ( n m ) = \gamma_3( n) m , \]
for any $ m \in \mathcal{M} $ and $ n \in \mathcal{N} $.

Set $ \lambda = \begin{bmatrix}
\bar{\alpha} ( 1_\mathcal{A} ) & 0 \\
0 & \bar{\beta} (1_\mathcal{B} ) 
\end{bmatrix} $. Since $ \eta ( \bar{\alpha} ( 1_\mathcal{A} ) ) = \bar{\beta} (1_\mathcal{B} ) $, it follows that $ \lambda \in Z(\mathcal{U} ) $. Now, we define linear map $ \chi : \mathcal{U} \to \mathcal{U} $ by 
\[ \chi \left( \begin{bmatrix}
a & m \\
n & b
\end{bmatrix} \right) = \begin{bmatrix}
\eta^{-1} ( \alpha_4(a)) + \beta_1 (b) & 0 \\
0 & \alpha_4(a) + \eta ( \beta_1 (b))
\end{bmatrix} \]
By the hypothesis $ (ii) $, we see that $ \chi (A) \in Z ( \mathcal{U} ) $ for any $ A \in \mathcal{U} $. 
Now, for any $ A= \begin{bmatrix}
a & m \\
n& b
\end{bmatrix} \in Z ( \mathcal{U} ) $, we can obtain
\begin{align*}
\phi (A) & = \phi \left( \begin{bmatrix}
a & m \\
n& b
\end{bmatrix} \right) = \begin{bmatrix}
\alpha_1 ( a) + \beta_1 (b) & \tau_2 (m) \\
\gamma_3 (n) & \alpha_4 (a) + \beta_4 (b)
\end{bmatrix} \\
& = \begin{bmatrix}
\bar{\alpha} (a) & \tau_2 (m) \\
\gamma_3 (n) & \bar{\beta} (b)
\end{bmatrix} + \begin{bmatrix}
\eta^{-1} ( \alpha_4(a)) + \beta_1 (b) & 0 \\
0 & \alpha_4(a) + \eta ( \beta_1 (b))
\end{bmatrix} \\
& = \begin{bmatrix}
\bar{\alpha}(1_\mathcal{A})a &\bar{\alpha}(1_\mathcal{A}) m \\
\bar{\beta} (1_\mathcal{B})n & \bar{\beta} (1_\mathcal{B} )b
\end{bmatrix} + \chi \left( \begin{bmatrix}
a & m \\
n & b
\end{bmatrix} \right) \\
& = \lambda A + \chi (A)
\end{align*}
Finally, by using $ \phi $ is Lie triple centralizer and the above results for any $ A, B, C \in \mathcal{U} $, we get
\begin{align*}
\chi ( [ [ A, B ] , C ] ) & = \phi ( [ [ A, B], C] ) - \lambda ( [ [ A, B], C] ) \\
& = ( [ [ \phi(A), B], C] )- \lambda ( [ [ A, B], C] ) \\
& = [ [ \lambda A + \chi (A) , B], C] ) - \lambda ( [ [ A, B], C] ) = 0 
\end{align*}
$ (iii) \Rightarrow (ii) $. It is clear.\\
$ (ii) \Rightarrow (iii) $. Let 
\[ a_0 = \alpha_1(1_\mathcal{A}) - \eta^{-1} ( \alpha_4 ( 1_\mathcal{A} ) ) , \]
\[ b_0 = \beta_4 (1_\mathcal{B} ) - \eta ( \beta _1 ( 1_\mathcal{B} )) . \]
By the assumption, and ($ iii $) and ($ iv $) of Theorem \ref{t1}, for any $ a \in \mathcal{A} $, any $ m \in \mathcal{M} $ and any $ n \in \mathcal{N} $, we have
\[ \tau_2(m) = a_0 m . \]
Therefore
\[ a \tau_2 (m) = a a_0 m = \alpha_1 (a) m - m \alpha_4 (a) . \]
So 
\[ ( \alpha_1(a) - a a_0 ) m = m \alpha_4 (a) , \]
and also $ \gamma_3 (n) = n a_0 $. Thus
\[ \gamma_3 (n) a = n a_0 a = n \alpha_1 (a) - \alpha_4 (a)n . \]
Hence
\[ n ( \alpha_1(a) - a_0 a) = \alpha_4 (a) n . \]
Now, by using Lemma \ref{l24}, we get
\[ \alpha_4 (a) \in \pi_\mathcal{B} ( Z(\mathcal{U})), ~~~ \text{and} ~~~ \eta^{-1} (\alpha_4(a)) = \alpha_1 (a) - a_0 a , \]
for any $ a \in \mathcal{A} $. 

By using the assumption, Theorem \ref{t1}-($ iii $) and Theorem \ref{t1}-($ iv $), for any $ b \in \mathcal{B} $, any $ m \in \mathcal{M} $ and any $ n \in \mathcal{N} $, we see that
\[ \tau_2(m) = m b_0 , \quad \gamma_3 (n) = b_0 n , \]
and
\[ \tau_2 (m) b = m b_0 b = m \beta_4 (b) - \beta_1 (b) m . \]
So 
\[ m ( \beta_4 (b) - b_0 b ) = \beta _1 (b) m , \]
also we have
\[ b \gamma_3 (n) = b b_0 n = \beta_4 (b) - n \beta_1(b) . \]
It follows that
\[ ( \beta_4 (b) - b b_0 ) n = n \beta_1 (b) . \]
By Lemma \ref{l24}, it implies $ \beta_1 (b) \in \pi_\mathcal{A} ( Z(\mathcal{U} )) $ and $ \eta ( \beta_1(b)) = \beta_4 
(b) - b b_0 $ for any $ b \in \mathcal{B} $. This completes the proof.
\end{proof}
In the following corollary, we provide sufficient conditions to characterize Lie triple centralizers.
\begin{cor}\label{c36}
Let $ \mathcal{U} $ satisfies 
\[ a \in \mathcal{A} , ~~ ~~a \mathcal{M} =0 ~~~ \text{and} ~~~ \mathcal{N} a = 0 \Rightarrow a = 0 , \]
\[ b \in \mathcal{B} , ~~~~ \mathcal{M} b =0 ~~~ \text{and} ~~~ b \mathcal{N} = 0 \Rightarrow b =0 . \]
Suppose that
\begin{enumerate}
\item[(i)] $\pi_\mathcal{B} ( Z(\mathcal{U} )) = Z(\mathcal{B} )$ or $[ [ \mathcal{A} , \mathcal{A} ] , \mathcal{A} ]= \mathcal{A}$;
\item[(ii)] $\pi_\mathcal{A} ( Z(\mathcal{U} )) = Z(\mathcal{A} )$ or $ [ [ \mathcal{B} , \mathcal{B} ], \mathcal{B} ] = \mathcal{B}$ .
\end{enumerate}
Then a linear mapping $ \phi : \mathcal{U} \to \mathcal{U} $ is a Lie triple centralizer if and only if $ \phi $ is a proper Lie triple centralizer.
\end{cor}
\begin{proof}
Let $ \phi : \mathcal{U} \to \mathcal{U} $ to be a Lie triple centralizer. If $ \pi_\mathcal{B} ( Z(\mathcal{U} )) = Z(\mathcal{B} ) $, then by using Corollary \ref{c2}, we see that $ \alpha_4 ( \mathcal{A} ) \subseteq Z ( \mathcal{B} ) $ and hence $ \alpha_4(\mathcal{A} ) \subseteq \pi_\mathcal{B} ( Z (\mathcal{U})) $. If $[ [ \mathcal{A} , \mathcal{A} ] , \mathcal{A} ]= \mathcal{A} $, by Corollary \ref{c2}, since $ \alpha_4 ( [ [ \mathcal{A} , \mathcal{A} ] , \mathcal{A} ] ) = 0 $, so $ \alpha_4 = 0 $ and it is clear that $ \alpha_4 ( \mathcal{A} ) \subseteq \pi_\mathcal{B} ( Z (\mathcal{U})) $. If one of the conditions (ii) is also holds, similarly from Corollary \ref{c2}, it follows that $ \beta_1 (\mathcal{B} ) \subseteq \pi_\mathcal{A} ( Z(\mathcal{U} )) $. Now we get the result from Theorem \ref{t2}. The converse is clear.
\end{proof}
\section{Some Applications}
In this section, we refer to some applications of the results obtained. First, as an application of Theorem \ref{t2}, we characterize the generalized Lie triple derivations on unital generalized matrix algebras. Throughout this section $ \mathcal{U} = \begin{bmatrix}
\mathcal{A} & \mathcal{M} \\
\mathcal{N} & \mathcal{B}
\end{bmatrix} $ is a unital 2-torsion free generalized matrix algebra. It is recalled that a linear map $d:\mathcal{U}\rightarrow \mathcal{U}$ is a \textit{Jordan derivation} if $d(A\circ B)=A\circ d(B)+ d(A)\circ B$ for every $A,B \in \mathcal{U}$. Also, $d$ is said to be a \textit{a singular Jordan derivation} if 
\[ d \left( \begin{bmatrix}
0 & m \\
n & 0
\end{bmatrix} \right) = \begin{bmatrix}
0 & \varrho(n) \\
\tau(m) & 0
\end{bmatrix}, \]
where $\tau: \mathcal{M}\rightarrow \mathcal{N}$ and $\varrho: \mathcal{N}\rightarrow \mathcal{M}$ are linear mappings (see \cite{ben}).
\par 
To prove the our result, we need the following theorem, which is proved in \cite{ben}. 
\begin{thm} \label{tah} $($\cite[Theorem 5.1]{ben}$)$ 
Let $\mathcal{U}$ satisfies 
\[ a \in \mathcal{A} , ~~ ~~a \mathcal{M} =0 ~~~ \text{and} ~~~ \mathcal{N} a = 0 \Rightarrow a = 0 , \]
\[ b \in \mathcal{B} , ~~~~ \mathcal{M} b =0 ~~~ \text{and} ~~~ b \mathcal{N} = 0 \Rightarrow b =0 . \]
Suppose that one of the following statements holds: 
\begin{enumerate}
\item[(i)] $[ [ \mathcal{A} , \mathcal{A} ] , \mathcal{A} ]= \mathcal{A}$ and $ [ [ \mathcal{B} , \mathcal{B} ], \mathcal{B} ] = \mathcal{B}$;
\item[(ii)] $\pi_\mathcal{A} ( Z(\mathcal{U} )) = Z(\mathcal{A} )$ and $[ [ \mathcal{A} , \mathcal{A} ] , \mathcal{A} ]= \mathcal{A}$;
\item[(iii)] $\pi_\mathcal{B} ( Z(\mathcal{U} )) = Z(\mathcal{B} )$ and $ [ [ \mathcal{B} , \mathcal{B} ], \mathcal{B} ] = \mathcal{B}$;
\item[(iv)] $\pi_\mathcal{A} ( Z(\mathcal{U} )) = Z(\mathcal{A} )$ and $\pi_\mathcal{B} ( Z(\mathcal{U} )) = Z(\mathcal{B} )$ and $\mathcal{A}$ or $\mathcal{B}$ satisfies 
\[ [x,\mathcal{G}]\in Z(\mathcal{G})\Rightarrow x \in Z(\mathcal{G})\]
for all $x\in \mathcal{G}$ where $\mathcal{G}\in \lbrace\mathcal{A}, \mathcal{B}\rbrace$.
\end{enumerate}
We also assume that one of the following statements holds: 
\begin{enumerate}
\item[(a)] $\mathcal{A}$ contains no central ideals;
\item[(b)] $\mathcal{B}$ contains no central ideals;
\item[(c)] $ Z ( \mathcal{U} ) = \left\lbrace \begin{bmatrix}
a & 0 \\
0 & b
\end{bmatrix} : a \in Z(\mathcal{A}), b \in Z(\mathcal{B}), am_0 =m_0 b \right\rbrace $ for some $m_0 \in \mathcal{M}$.
\item[(d)] $ Z ( \mathcal{U} ) = \left\lbrace \begin{bmatrix}
a & 0 \\
0 & b
\end{bmatrix} : a \in Z(\mathcal{A}), b \in Z(\mathcal{B}), n_0 a=bn_0 \right\rbrace $ for some $n_0 \in \mathcal{N}$.
\end{enumerate}
Then any Lie triple derivation $\xi:\mathcal{U} \rightarrow \mathcal{U}$ is of the form $\xi =\delta+ d +\Gamma$, where $\delta:\mathcal{U} \rightarrow \mathcal{U}$ is a derivation, $d:\mathcal{U} \rightarrow \mathcal{U}$ is a singular Jordan derivation and $\Gamma:\mathcal{U} \rightarrow Z(\mathcal{U})$ is a linear map that vanishes on $[[\mathcal{U},\mathcal{U}], \mathcal{U}]$. 
\end{thm}
The following theorem is the main result of this section, which is also a partial generalization of \ref{tah}. 
\begin{thm} \label{tg}
Let $\mathcal{U}$ satisfies 
\[ a \in \mathcal{A} , ~~ ~~a \mathcal{M} =0 ~~~ \text{and} ~~~ \mathcal{N} a = 0 \Rightarrow a = 0 , \]
\[ b \in \mathcal{B} , ~~~~ \mathcal{M} b =0 ~~~ \text{and} ~~~ b \mathcal{N} = 0 \Rightarrow b =0 . \]
Suppose that one of the following statements holds: 
\begin{enumerate}
\item[(i)] $[ [ \mathcal{A} , \mathcal{A} ] , \mathcal{A} ]= \mathcal{A}$ and $ [ [ \mathcal{B} , \mathcal{B} ], \mathcal{B} ] = \mathcal{B}$;
\item[(ii)] $\pi_\mathcal{A} ( Z(\mathcal{U} )) = Z(\mathcal{A} )$ and $[ [ \mathcal{A} , \mathcal{A} ] , \mathcal{A} ]= \mathcal{A}$;
\item[(iii)] $\pi_\mathcal{B} ( Z(\mathcal{U} )) = Z(\mathcal{B} )$ and $ [ [ \mathcal{B} , \mathcal{B} ], \mathcal{B} ] = \mathcal{B}$;
\item[(iv)] $\pi_\mathcal{A} ( Z(\mathcal{U} )) = Z(\mathcal{A} )$ and $\pi_\mathcal{B} ( Z(\mathcal{U} )) = Z(\mathcal{B} )$ and $\mathcal{A}$ or $\mathcal{B}$ satisfies 
\[ [x,\mathcal{G}]\in Z(\mathcal{G})\Rightarrow x \in Z(\mathcal{G})\]
for all $x\in \mathcal{G}$ where $\mathcal{G}\in \lbrace\mathcal{A}, \mathcal{B}\rbrace$.
\end{enumerate}
We also assume that one of the following statements holds: 
\begin{enumerate}
\item[(a)] $\mathcal{A}$ contains no central ideals;
\item[(b)] $\mathcal{B}$ contains no central ideals;
\item[(c)] $ Z ( \mathcal{U} ) = \left\lbrace \begin{bmatrix}
a & 0 \\
0 & b
\end{bmatrix} : a \in Z(\mathcal{A}), b \in Z(\mathcal{B}), am_0 =m_0 b \right\rbrace $ for some $m_0 \in \mathcal{M}$.
\item[(d)] $ Z ( \mathcal{U} ) = \left\lbrace \begin{bmatrix}
a & 0 \\
0 & b
\end{bmatrix} : a \in Z(\mathcal{A}), b \in Z(\mathcal{B}), n_0 a=bn_0 \right\rbrace $ for some $n_0 \in \mathcal{N}$.
\end{enumerate}
Let $ \Lambda: \mathcal{U} \to \mathcal{U} $ be a generalized Lie triple derivation associated with the Lie triple derivation $\xi: \mathcal{U} \to \mathcal{U}$. Then $ \Lambda (A) = \delta ( A ) +d(A) + \psi (A) + \lambda A $ for any $ A \in \mathcal{U} $, where $ \delta $ is a derivation on $ \mathcal{U} $, $d$ is a singular Jordan derivation on $\mathcal{U}$, $ \lambda \in Z(\mathcal{U} ) $ and $ \psi:\mathcal{U} \rightarrow Z(\mathcal{U})$ is a linear map that vanishes on $[[\mathcal{U},\mathcal{U}], \mathcal{U}]$.
\end{thm}
\begin{proof}
By Remark \ref{RE} the linear map $ \phi = \Lambda - \xi $ is a Lie triple centralizer on $ \mathcal{U} $. Thus by Corollary \ref{c36}, $ \phi (A) = \lambda A + \chi (A) $ for any $ A \in \mathcal{U} $, where $ \lambda \in Z (\mathcal{U} ) $ and $ \chi: \mathcal{U} \rightarrow Z (\mathcal{U} )$ is a linear mapping whith $ \chi ( [ [ A , B ] , C ] ) = 0 $ for all $ A , B , C \in \mathcal{U} $. By Theorem \ref{tah}, $ \xi = \delta+ d + \Gamma $, where $\delta$ is a derivation, $d$ is a singular Jordan derivation and $\Gamma:\mathcal{U} \rightarrow Z(\mathcal{U})$ is a linear map that $ \gamma ( [ [ A , B ] , C ] ) = 0 $ for all $ A , B , C \in \mathcal{U} $.
\par 
Suppose that $ \psi = \chi + \Gamma $. Thus $ \psi : \mathcal{U} \to Z(\mathcal{U}) $ is a linear map that $ \psi ( [ [ A , B ] , C ] ) = 0 $ for all $ A , B , C \in \mathcal{U} $, and 
\begin{align*}
\Lambda(A) & = \phi (A) + \Delta (A) \\
& = \lambda A + \chi (A) + \delta (A) +d(A)+ \Gamma (A) \\
&= \delta (A)+d(A)+ \psi(A) + \lambda A
\end{align*}
for any $ A \in \mathcal{U} $. The proof is complete. 
\end{proof}
In the following remarks, we refer to applications of the results obtained to some special algebras.
\begin{rem}\label{tri}
Let $ \mathcal{T} = T ri(\mathcal{A}, \mathcal{M}, \mathcal{B})$ be a unital triangular algebra. Suppose that $ \mathcal{M} $ is a faithful ($\mathcal{A}, \mathcal{B}$)-bimodule. Then the results obtained in Theorem \ref{t2}, Corollary \ref{c36} and Theorem \ref{tg} are valid for $ \mathcal{T}$.
\end{rem}
\begin{rem}\label{oa}
The conditions of Theorem \ref{t2}, Corollary \ref{c36} and Theorem \ref{tg} hold for each of the following algebras. The definition of these algebras and the proof that the conditions of Theorem \ref{t2}, Corollary \ref{c36} and Theorem \ref{tg} are true for these algebras can be found in the mentioned references.
\begin{enumerate}
\item[(i)] $M_{n}(\mathcal{A}) $ ($n \geq 2$) the algebra of all $n\times n$ matrices over a unital algebra $\mathcal{A}$ (see \cite{xi});
\item[(ii)] $ T_n (\mathcal{A} ) $ ($n \geq 2$) the algebra of all $ n \times n $ upper triangular matrices with entries from a unital algebra $ \mathcal{A} $ (see \cite{xi});
\item[(iii)] Standard operator algebra $\mathcal{U}$ on a complex Banach space $\mathcal{X}$ (see \cite{qi1});
\item[(v)] Factor von Neumann algebra $\mathcal{M}$ on a complex Hilbert space $ \mathcal{H} $ (see \cite{qi1});
\item[(vi)] Non-trivial nest algebra $Alg\mathcal{N}$ on a complex Hilbert space $ \mathcal{H} $ (see \cite{xi}).
\end{enumerate}
\end{rem}
\begin{rem}
\begin{enumerate}
\item[(i)] It should be noted that using the method of proving Theorem \ref{tah}, the results obtained for Lie triple derivations (for example, the results obtained in \cite{ben}) can be generalized to generalized Lie triple derivations on generalized matrix algebras.
\item[(ii)] Since any Lie centralizer and Jordan centralizer is a Lie triple centralizer, by using the results obtained in this paper, Lie centralizers and Jordan centralizers can also be characterized on generalized matrix algebras. Therefore, with the method used to prove Theorem \ref{tah}, the results obtained for Lie derivations can be generalized to generalized Lie derivations on generalized matrix algebras (See before Remark \ref{RE}).
\end{enumerate}
\end{rem}

\subsection*{Acknowledgment}
The authors thanks the referees for careful reading of the manuscript and for helpful suggestions.

\bibliographystyle{amsplain}
\bibliography{xbib}

\end{document}